\documentclass[a4paper,11pt]{article}
\title{Isometric embeddings of a class of separable metric spaces into Banach spaces}
\author{S.K.Mercourakis and G.Vassiliadis}
\date{}
\usepackage{amsmath}
\usepackage{amsfonts}
\usepackage{amssymb}
\usepackage{amscd}
\usepackage{amsthm}
\usepackage{makeidx}
\usepackage{euscript}
\theoremstyle{plain}
\newtheorem{theo}{Theorem}
\newtheorem{lemm}{Lemma}
\newtheorem{prop}{Proposition}

\theoremstyle{definition}

\newtheorem*{claim}{\underline{Claim}}
\newcommand{\beg}{\begin{proof}[Proof]}
\newcommand{\eq}{equilateral $\,$}

\newcommand{\lap}{\ell_\infty}

\begin{document}
\maketitle

\begin{abstract}
\footnotesize 
Let $(M,d)$ be a bounded countable metric space and $c>0$ a constant, such that
$d(x,y)+d(y,z)-d(x,z) \ge c$, for any pairwise distinct points $x,y,z$ of $M$.
For such metric spaces we prove that they can be isometrically embedded into
any Banach space containing an isomorphic copy of $\lap$.

\normalsize 
\end{abstract}


\section*{Introduction}
Let $(M,d)$ be a metric space; following \cite{K} we will call it \textit{concave}, when 
the triangle inequality is strict, i.e. when $d(x,y)+d(y,z)>d(x,z)$ for any pairwise 
distinct points $x,y,z$ of $M$.

In this note we are interested in (concave) metric spaces satisfying the stronger property:
there is a constant $c>0$, such that
$d(x,y)+d(y,z)-d(x,z) \ge c$, for any pairwise distinct points $x,y,z$. Let us call
these spaces \textit{strongly concave} metric spaces.

The main result we prove is an infinite dimensional version of Theorem 4.3 of \cite{K}, 
that is, if a Banach space $X$ contains an isomorphic copy of $\lap$, then $X$ contains
isometrically any bounded countable strongly concave metric space (Th.2). An immediate
consequence of this result is that any Banach space containing an isomorphic copy of $c_0$,
admits an infinite equilateral set (Th.3). This result was first proved (by similar methods)
in \cite{MV} (Th.2).

A subset $S$ of a metric space $(M,d)$ is said to be equilateral, if there is a $\lambda>0$
such that for $x \neq y \in S$ we have $d(x,y)=\lambda$; we also call $S$ a $\lambda$-equilateral
set (see \cite{S}).

If $X$ is any (real) Banach space, then $B_X$ and $S_X$ denote its closed unit ball
and unit sphere respectively. $X$ is said to be strictly convex, if for any $x \neq y \in S_X$
we have $\|x+y\|<2$. The Banach-Mazur distance between two isomorphic Banach spaces $X$ and $Y$
is $d(X,Y)=\inf \{\|T\| \cdot \|T^{-1}\|: T \:\text{is an isomorphism} \}$.

\section*{Strongly concave metric spaces}

We start by presenting some examples of concave metric spaces

\noindent \textbf{Examples 1} \\
\noindent (1) a) Let $(M,d)$ be a discrete metric space (i.e. $d(x,y)=1$ when $x \neq y$).
Clearly $1=d(x,z)<d(x,y)+d(y,z)=2$ for any pairwise distinct triplet $x,y,z \in M$. 
Therefore $(M,d)$ is a concave metric space.
In particular, every $\lambda$-equilateral subset of any metric space is a concave 
metric space.\\
b) More generally, every \textit{ultrametric} space is concave. This holds since
for any $x,y,z$ pairwise distinct points we have $d(x,z) \le \max \{d(x,y),d(y,z)\} <
d(x,y)+d(y,z)$.

\noindent (2) Let $(X,\| \cdot \|)$ be a strictly convex Banach space. As is well known, if 
$x,y,z$ are non collinear points of $X$ then $\|x-z\|<\|x-y\|+\|y-z\|$.\\
It then follows that the unit sphere $S_X$ and every affinely independent subset $A$ of $X$
with the norm metric are concave metric spaces (in any case no three pairwise distinct points 
are collinear). 

\noindent (3) Let $(X,\| \cdot \|)$ be a Banach space and $A \subseteq B_X$ such that 
$x \neq y \in A \Rightarrow \|x-y\|>1$ (see \cite{GM}). Then for any $x,y,z$ pairwise distinct
points of $A$ we have $\|x-y\|+\|y-z\|-\|x-z\|>1+1-\|x-z\| \ge 1+1-2=0$.
Hence $A$ with the norm metric is concave.

\noindent (4)  Let $(M,d)$ be any metric space and $p \in (0,1)$. Then it is
rather easy to show that $d^p$ is a concave metric on $M$. This follows from
the fact that, given $a,b,c>0$ with $a \le b+c$ then $a^p<b^p+c^p$.
The metric $d^p$ is then called the snowflaked version of $d$ (see \cite{O}).

We are interested in concave metric spaces $(M,d)$ satisfying the stronger property:
there is a constant $c>0$ such that for any pairwise distinct points $x,y,z$
of $M$ we have $d(x,y)+d(y,z)-d(x,z) \ge c$, equivalently $d(x,z)+c \le d(x,y)+d(y,z)$.
Let us call these spaces \textit{strongly concave} spaces.

\begin{lemm}
Every strongly concave metric space is separated (or uniformly discrete).
\end{lemm}

\beg 
Assume that $(M,d)$ is a $c$-strongly concave metric space.
We claim that $x \neq y \in M \Rightarrow d(x,y) \ge \frac{c}{2}$.
Assume for the purpose of contradiction that there is a pair $\{x,y\}
\subseteq M$ with $d(x,y)<\frac{c}{2}$. Let also $z \in M \setminus \{x,y\}$.
We then have $d(x,y)+d(y,z) \le d(x,y)+(d(y,x)+d(x,z))=2 d(x,y)+d(x,z) \Rightarrow$
$d(x,y)+d(y,z)-d(x,z) \le 2 d(x,y)<2 \frac{c}{2}=c$. The last inequality clearly
contradicts the fact that $M$ is $c$-strongly concave.

\end{proof}

The following are examples of strongly concave metric spaces.

\noindent \textbf{Examples 2} \\
\noindent (1) Every finite concave metric space is clearly strongly concave.

\noindent (2) Let $A$ be a a $\lambda$-\eq subset of any metric space $(M,d)$.
For any pairwise distinct points $x,y,z$
of $A$ we have $d(x,y)+d(y,z)-d(x,z)=\lambda+\lambda-\lambda=\lambda$, so $A$ 
is a $\lambda$-strongly concave metric subspace of $(M,d)$.

\noindent (3) Let $(X,\| \cdot \|)$ be a Banach space. Also let $A \subseteq B_X$
with the property that $x \neq y \in A \Rightarrow \|x-y\| \ge 1+\varepsilon$,
where $\varepsilon>0$ is a constant. Then we have $\|x-y\|+\|y-z\|-\|x-z\|>(1+\varepsilon)+
(1+\varepsilon)-2=2 \varepsilon$ (cf. Examples 1(3)). Therefore $A$ with the norm metric 
is a $2 \varepsilon$-strongly concave metric space.

Note that if $\dim X=\infty$, then by a result of Elton and Odell (\cite{EO}) there is 
$A \subseteq S_X$ infinite and $\varepsilon>0$ such that 
$x \neq y \in A \Rightarrow \|x-y\| \ge 1+\varepsilon$.

\noindent \textbf{Remarks 1} 
\noindent (1) Clearly every separable strongly concave metric space $M$ is at most countable (this is so
because $M$ is separated, hence it has the discrete topology).

\noindent (2) Every subspace of a concave (resp. strongly concave) space has the same property.

The following result is classical (see \cite{O}).

\begin{theo} (Frech\`{e}t) Every  separable metric space $(M,d)$ embeds isometrically
into $\ell_\infty$.

\end{theo}

\beg
Let $(x_n) \subseteq M$ be a dense sequence in $M$. Then the map
\[\varphi:x \in M \mapsto (d(x,x_n)-d(x_1,x_n))_{n \ge 1} \in \ell_\infty \]
satifies our claim.

\end{proof}

\noindent \textbf{Remark 2} Let $(M,d)$ be a separable metric space. We define a map
\[\sigma:M \rightarrow \Bbb{R}^\Bbb{N} \;\:\text{with}\:\; \sigma(x)=(d(x,x_n))_{n \ge 1} \]
where $(x_n)$ is any dense sequence in $M$. Then the Frech\`{e}t embedding
of $M$ into $\ell_\infty$ is the map
\[\varphi(x)=\sigma(x)-\sigma(x_1),\, x\in X\]
Note that if the space $(M,d)$ is bounded (that is, there is $k>0$ such that $d(x,y) \le k$
for all $x,y \in M$), then the map $\sigma$ is already an isometric embedding of $M$ into $\lap$,
which we will still call the Frech\`{e}t embedding of $M$ into $\ell_\infty$.

\begin{prop}
Let $(M,d)$ be a bounded countable infinite metric space. Then there is an infinite
subset $N$ of $M$ such that the Frech\`{e}t embedding of $N$ into $\ell_\infty$
takes values into the space $\mathbf{c}$.
\end{prop}

\beg 
Let $\{x_1,x_2,\dots,x_n,\dots\}$ be a one-to-one enumeration of $M$. Then
$\sigma(x_k)=(d(x_k,x_n))_{n \ge 1} \in \ell_\infty$, for $k \in \Bbb{N}$, 
since $d$ is a bounded metric. We construct by induction a subsequence
$\{x'_1,x'_2,\dots,x'_n,\dots\}$ of $(x_n)$ satisfying our claim. 

Since $(d(x_1,x_n))_{n \ge 1}$ is a bounded sequence of real numbers, there is
$A_1 \subseteq \mathbb{N}$ infinite, such that $d(x_1,x_n) \stackrel{n \in A_1}{\longrightarrow} \alpha_1$.
Set $n_1=1$.

Let $n_2=\min A_1$, for which we may assume that $n_2>n_1$. Then for the sequence $(d(x_{n_2},x_n))_{n \in A_1}$,
there is $A_2 \subseteq A_1$ infinite with $n_3 =\min A_2>n_2$ such that 
$d(x_{n_2},x_n) \stackrel{n \in A_2}{\longrightarrow} \alpha_2$.

Then for the sequence $(d(x_{n_3},x_n))_{n \in A_2}$,
there is $A_3 \subseteq A_2$ infinite with $n_4 =\min A_3>n_3$ such that 
$d(x_{n_3},x_n) \stackrel{n \in A_3}{\longrightarrow} \alpha_3$.

The inductive process should be clear. Now set $A=\{n_1<n_2<\dots<n_k<\dots\}$.
Clearly $\{n_k,n_{k+1},\dots \} \subseteq A_k$ for $k \ge 1$ and hence 
$d(x_{n_k},x_n) \stackrel{n \in A} \longrightarrow \alpha_k$ for all $k \ge 1$.
It is clear that the set $N=\{x'_k=x_{n_k}:k \ge 1\}$ satisfies our requirements.
\end{proof}

The following theorem is the main result of this note; its proof resembles the
proof of Theorem 4.3 of \cite{K} and the proof of Theorem 2 of \cite{MV}
(we use Schauder's fixed point theorem the same way we did in \cite{MV}).
The origins of these ideas can be traced in Brass (see \cite{B} and \cite{S}) and
Swanepoel and Villa (see \cite{SV1} and \cite{SV2}).

\begin{theo}
Let $X$ be any Banach space containing an isomorphic copy of $\ell_\infty$. Then $X$
contains isometrically any bounded separable strongly concave metric space.
\end{theo}

\beg
We shall use a kind of non distortion property of $\ell_\infty$ proved
independently by Talagrand (\cite{T}) and Partington (\cite{Pa}). Let
us denote by $\|\cdot\|_\infty$ the usual norm of $\ell_\infty$.

\begin{claim}
Let $(M,d)$ be any bounded separable strongly concave metric space. There is 
$\delta>0$, such that if $\|\cdot\|$ is any equivalent norm on $\ell_\infty$ with
Banach Mazur distance \[d\left((\lap,\|\cdot\|_\infty),(\lap,\|\cdot\|)\right) \le 1+\delta \]
then the space $(M,d)$ embeds isometrically into $(\lap,\|\cdot\|)$.
\end{claim}

\underline{Proof of the Claim:}
Since $(M,d)$ is strongly concave, there is $\eta>0$ such that
$d(x,y)+d(y,z)-d(x,z) \ge \eta$, for each triplet $x,y,z$ of 
pairwise distinct points of $M$. We may assume that $\|x\| \le \|x\|_\infty
\le (1+\delta)\|x\|$ for $x \in \lap$, where $\delta>0$ is to be determined.

Let $I=\{(m,n):n<m,\;n,m \in \Bbb{N}\}$; denote by $K$ the compact cube $[0,\eta]^I$.
Since $M$ is (strongly concave and) separable, it is at most countable, so let
$M=\{x_1,x_2,\dots,x_n,\dots\}$. For $\varepsilon=(\varepsilon_{(m,n)}) \in K$ set
\[p_1(\varepsilon)=\left(d(x_1,x_1)-d(x_1,x_1),d(x_1,x_2)-d(x_1,x_2),\dots,d(x_1,x_n)-d(x_1,x_n),\dots \right)=(0,\dots,0,\dots)\]
\[p_2(\varepsilon)=\left(d(x_2,x_1)-d(x_1,x_1)+\varepsilon_{(2,1)},d(x_2,x_2)-d(x_1,x_2),\dots,d(x_2,x_n)-d(x_1,x_n),\dots \right)\]
\vdots
\[p_n(\varepsilon)=\left(d(x_n,x_1)-d(x_1,x_1)+\varepsilon_{(n,1)},\dots,d(x_n,x_{n-1})-d(x_1,x_{n-1})
+\varepsilon_{(n,n-1)},d(x_n,x_n)-d(x_1,x_n),\dots \right)\]
\vdots

\noindent (Note that $x_n \mapsto p_n(0)$ is the Frech\`{e}t embedding of $M$ into $(\lap,\|\cdot\|_\infty)$).

For $n<m$ we have 
\[\|p_n(\varepsilon)-p_m(\varepsilon)\|_\infty= \sup_k |d(x_n,x_k)+\varepsilon_{(n,k)}-
(d(x_m,x_k)+\varepsilon_{(m,k)})|\] 
where we set $\varepsilon_{(k,l)}=0$, for $l \ge k$. This supremum is equal
to $d(x_n,x_m)+\varepsilon_{(m,n)}$, as for $k \neq n,m$ we have
\[d(x_n,x_k)-d(x_m,x_k)+\varepsilon_{(n,k)}-\varepsilon_{(m,k)} \le d(x_n,x_m)- \eta+\varepsilon_{(n,k)}-
\varepsilon_{(m,k)} \le d(x_n,x_m).\]

We define a function
\[\varepsilon=\left(\varepsilon_{(m,n)}\right) \in K \stackrel{\varphi}{\mapsto}
\varphi(\varepsilon)=\left(\varphi_{(m,n)}(\varepsilon)\right) \in K,\]
by the rule $\varphi_{(m,n)}(\varepsilon)=d(x_n,x_m)+\varepsilon_{(m,n)}-
\|p_n(\varepsilon)-p_m(\varepsilon)\|$. Note that
$\varphi_{(m,n)}(\varepsilon) \ge d(x_n,x_m)+\varepsilon_{(m,n)}-
\|p_n(\varepsilon)-p_m(\varepsilon)\|_\infty=0$ (using the computation above and the fact
that the norm $\|\cdot\|_\infty$ dominates $\|\cdot\|$).
We also have
\[d(x_n,x_m)+\varepsilon_{(m,n)}=\|p_n(\varepsilon)-p_m(\varepsilon)\|_\infty \le
(1+\delta) \|p_n(\varepsilon)-p_m(\varepsilon)\| \Rightarrow \]
\[ \frac{1}{1+\delta} (d(x_n,x_m)+\varepsilon_{(m,n)}) \le \|p_n(\varepsilon)-p_m(\varepsilon)\|. \]
Therefore \[\varphi_{(m,n)}(\varepsilon)=d(x_n,x_m)+\varepsilon_{(m,n)}-
\|p_n(\varepsilon)-p_m(\varepsilon)\| \]
\[\le d(x_n,x_m)+\varepsilon_{(m,n)}-\frac{1}{1+\delta} \left(d(x_n,x_m)+\varepsilon_{(m,n)}\right)  \]
\[=\frac{\delta}{1+\delta} \left(d(x_n,x_m)+\varepsilon_{(m,n)}\right). \]
It then follows from (this inequality and) the fact that $M$ is bounded that if $\delta$
is quite small, then $\varphi_{(m,n)}(\varepsilon) \le \eta$, for $\varepsilon \in K$.

Since each coordinate function $\varphi_{(m,n)}$ is continuous (as dependent on finite coordinates, 
i.e. from the set $\{(k,l):1 \le l<k \le m\}$) it follows that $\varphi$ is also continuous. By a classical
result of Schauder, $\varphi$ has a fixed point $\varepsilon'=\left(\varepsilon'_{(m,n)}\right) \in K$,
that is $\varphi(\varepsilon')=\varepsilon'$, which implies $\|p_n(\varepsilon')-p_m(\varepsilon')\|=d(x_n,x_m)$,
for all $n,m \in \Bbb{N}$. The proof of the Claim is complete. 

Denote by $\|\cdot\|$ the norm of $X$ and let $Y$ be a subspace of $X$ isomorphic to $\lap$.
By the non distortion property of $(\lap,\|\cdot\|_\infty)$ there is a subspace $Z \subseteq Y$
(isomorphic to $\lap$) such that \[d\left((Z,\|\cdot\|),(\lap,\|\cdot\|_\infty)\right) \le 1+\delta\]
(this is the $\delta>0$ postulated in the Claim).
It follows immediately from the Claim that the space $(Z,\|\cdot\|)$ contains an isometric copy of $(M,d)$.

\end{proof}

In the special case when $(M,d)$ is the countable infinite discrete metric space we get the following
result first proved in \cite{MV} (Th.2), essentially with the same method 

\begin{theo}
Every Banach space $X$ containing an isomorphic copy of $c_0$ admits an infinite equilateral set.
\end{theo}

\beg
Take in the proof of the previous theorem $(M,d)$ to be the countable infinite discrete space.
Then $\eta=1$ and the resulting family $(p_n(\varepsilon))_{n \ge 1},\;\varepsilon \in K=
\left[0,1 \right]^I$ takes values in $c_0$ (remember that $x_n \mapsto p_n(0)$ is the 
Frech\`{e}t embedding of $(M,d)$ into $c_0$). Since $(c_0,\|\cdot\|_\infty)$ is non distortable, we get 
the conclusion.
\end{proof}

Theorem 2 can be improved in the following way

\begin{theo}
Let $(M,d)$ be an infinite bounded separable strongly concave metric space. Then there is $N \subseteq M$
infinite such that the metric space $(N,d)$ can be isometrically embedded into any Banach space containing
an isomorphic copy of the space $c_0$.
\end{theo}

\beg
By Proposition 1, there is $N \subseteq M$ infinite such that the Frech\`{e}t embedding $\sigma:N \rightarrow \lap$
takes values into $\mathbf{c}$. Then the proof of Theorem 2 gives us a family of embeddings 
$(p_n(\varepsilon))_{n \ge 1},\;\varepsilon \in K=[0,\eta]^I$ taking values into $\mathbf{c}$.
Since $\mathbf{c}$ is isomorphic to $c_0$, we are done.
\end{proof}

\noindent S.K.Mercourakis, G.Vassiliadis\\
University of Athens\\
Department of Mathematics\\
15784 Athens, Greece\\
e-mail: smercour@math.uoa.gr

\hspace{0.7cm} georgevassil@hotmail.com
\end{document}